\documentclass[11pt, a4paper]{article}
\usepackage{amsmath}
\usepackage{amsfonts}
\usepackage{amssymb}
\usepackage{ulem}

\usepackage{amsthm}
\usepackage{url}

\usepackage{cite}
\usepackage{fancyhdr}

\usepackage{graphicx, color}
\usepackage{hyperref}

\newtheorem{thm}{Theorem}

\newtheorem{prop}{Proposition}

\newtheorem*{rem*}{Remark}
\newtheorem{sta}{Statement}

\begin{document}

\begin{center}

\textsf{ \Large{Properties of Convex Pentagonal Tiles for Periodic Tiling}}

\vspace{0.5cm}

Teruhisa S{\small UGIMOTO}

\vspace{0.3cm}

{\footnotesize

\textit{The Interdisciplinary Institute of Science, Technology and Art,\\
Suzukidaini-building 211, 2-5-28 Kitahara, Asaka-shi, Saitama, 351-0036, 
Japan}

\textit{E-mail address:}
\textit{ismsugi@gmail.com}

}
\end{center}

{\small
\begin{abstract}
\noindent
A convex pentagonal tile is a convex pentagon that 
admits a monohedral tiling. We show that a convex pentagonal tile that 
admits a periodic tiling has a property in which the sum of three 
internal angles of the pentagon is equal to $360^ \circ$.
\end{abstract}
}


\section{Introduction}
\label{section1}

A \textit{tiling} or \textit{tessellation} of the plane is a collection of sets 
that are called tiles, which covers a plane without gaps and overlaps, 
except for the boundaries of the tiles. If all the tiles in a tiling are of 
the same size and shape, then the tiling is \textit{monohedral}. Then, a 
polygon that admits a monohedral tiling is called a \textit{prototile} of 
monohedral tiling, or simply, a\textit{ polygonal tile}. A tiling by convex 
polygons is \textit{edge-to-edge} if any two convex polygons in a tiling 
are either disjoint or share one vertex or an entire edge in 
common~\cite{G_and_S_1987, Hallard_1991, Sugimoto_2012a, Sugimoto_2015, 
Sugimoto_NoteTP, Sugimoto_2016, Sugimoto_2017rcmm1, Sugimoto_2017rcmm2}.
A tiling of the plane is \textit{periodic} if the tiling can be translated 
onto itself in two nonparallel directions~\cite{G_and_S_1987, Hallard_1991, 
Sugimoto_2017rcmm1}.

\begin{figure}[htbp]
 \centering\includegraphics[width=15.5cm,clip]{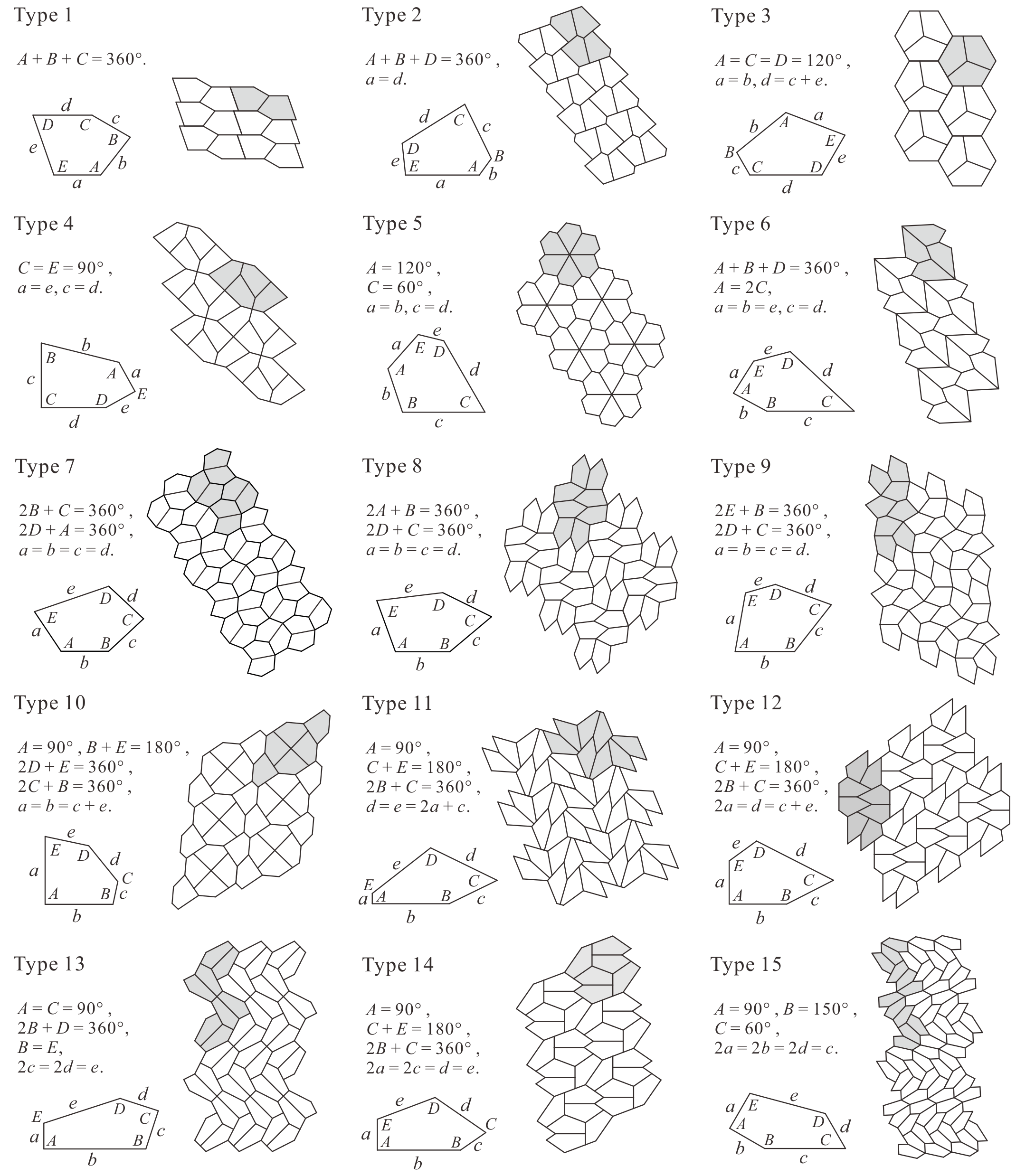} 
  \caption{{\small 
Convex pentagonal tiles of 15 families. Each of the convex 
pentagonal tiles is defined by some conditions between the lengths of the 
edges and the magnitudes of the angles, but some degrees of freedom remain. 
For example, a convex pentagonal tile belonging to Type 1 satisfies that the 
sum of three consecutive angles is equal to $360^ \circ$. This condition for 
Type 1 is expressed as $A+B+C=360^ \circ$ in this figure. The pentagonal tiles of 
Types 14 and 15 have one degree of freedom, that of size. For example, the 
value of $C$ of the pentagonal tile of Type 14 is $\cos ^{ - 1}((3\sqrt {57} - 
17) / 16) \approx 1.2099\;$rad $ \approx 69.32^ \circ $. The pale gray 
pentagons in each tiling indicate a fundamental region (the unit that can 
generate a periodic tiling by translation only).
}
\label{fig1}
}
\end{figure}

In the classification problem of convex polygonal tiles, only the pentagonal 
case is open. At present, fifteen families of convex pentagonal tiles, each 
of them referred to as a ``Type'' (i.e., Type 1, Type 2, etc. up to Type 
15), are known (see Figure~\ref{fig1}) but it is not known whether this list is 
complete\footnote{ In May 2017, Micha\"{e}l Rao declared that the complete 
list of Types of convex pentagonal tiles had been obtained (i.e., they have 
only the known 15 families), but it does not seem to be fixed as of October 
2018~\cite{Rao_2017, wiki_pentagon_tiling }}~\cite{Gardner_1975a, Gardner_1975b, 
G_and_S_1987, Hallard_1991, Hirshh_1985, Kershner_1968, Klamkin_1980, Mann_2015, 
Rao_2017, Reinhardt_1918, Schatt_1978, Stein_1985, Sugimoto_2012a, 
Sugimoto_NoteTP,Sugimoto_2017rcmm2, Wells_1991, wiki_pentagon_tiling}. 
However, it has been proved that a convex pentagonal tile that can generate 
an edge-to-edge tiling belongs to at least one of the eight known 
types~\cite{Bagina_2011, Sugimoto_2012a, Sugimoto_2012b, Sugimoto_2015,
Sugimoto_NoteTP, Sugimoto_2016}. We are interested in the 
problem of convex pentagonal tiling (i.e., the complete list of Types of 
convex pentagonal tiles\footnote{ The classification of Types of 
convex pentagonal tiles is based on the essentially different properties of 
pentagons. The conditions of each Type express the essential properties. The 
classification problem of Types of convex pentagonal tiles and the 
classification problem of pentagonal tilings are quite different. The Types 
are not necessarily ``disjoint,'' that is, convex pentagonal tiles belonging 
to some Types also exist~\cite{Sugimoto_2012a, Sugimoto_2016}.}, regardless 
of edge-to-edge and non-edge-to-edge tilings). In this paper, we prove the following.

\begin{thm}\label{thm1}
If a convex pentagon is a convex pentagonal tile that 
admits a periodic tiling, the pentagon has a property in which 
the sum of three internal angles of the pentagon is equal to $360^ \circ$. 
\end{thm}

By using the notations $A$, $B$, $C$, $D$, and $E$ for vertices of convex pentagon, 
the combinations in which the sum of three internal angles is equal to  $360^ \circ$ 
can be expressed using the following equations:

\bigskip
\noindent
$A+B+C =360^ \circ $, $B+C+D = 360^ \circ $, $C+D+E =360^ \circ $, 
$D+E+A =360^ \circ $, $E+A+B = 360^ \circ $, $A+B+D = 360^ \circ $, 
$B+C+E =360^ \circ $, $C+D+A =360^ \circ $, $D+E+B =360^ \circ $, 
$E+A+C =360^ \circ $, $2A+B =360^ \circ $, $2A+C =360^ \circ $, 
$2A+D =360^ \circ $, $2A+E =360^\circ$, $2B+A =360^ \circ $, 
$2B+C =360^ \circ $, $2B+D =360^ \circ $, $2B+E =360^ \circ $, 
$2C+A =360^ \circ $, $2C+B =360^\circ$, $2C+D =360^ \circ $, 
$2C+E =360^\circ$, $2D+A =360^ \circ $, $2D+B = 360^ \circ $, 
$2D+C = 360^ \circ $, $2D+E = 360^ \circ $, $2E+A = 360^ \circ $, 
$2E+B = 360^ \circ $, $2E+C = 360^ \circ $, $2E+D =360^ \circ $, 
$3A =360^ \circ $, $3B =360^\circ$, $3C =360^ \circ $, 
$3D =360^\circ$, $3E = 360^ \circ $.
\bigskip

\noindent
That is, a convex pentagonal tile that admits a periodic tiling has at 
least one of the above relationships.

\renewcommand{\figurename}{{\small Figure.}}

\renewcommand{\figurename}{{\small Figure.}}
\begin{figure}[htbp]
 \centering\includegraphics[width=10cm,clip]{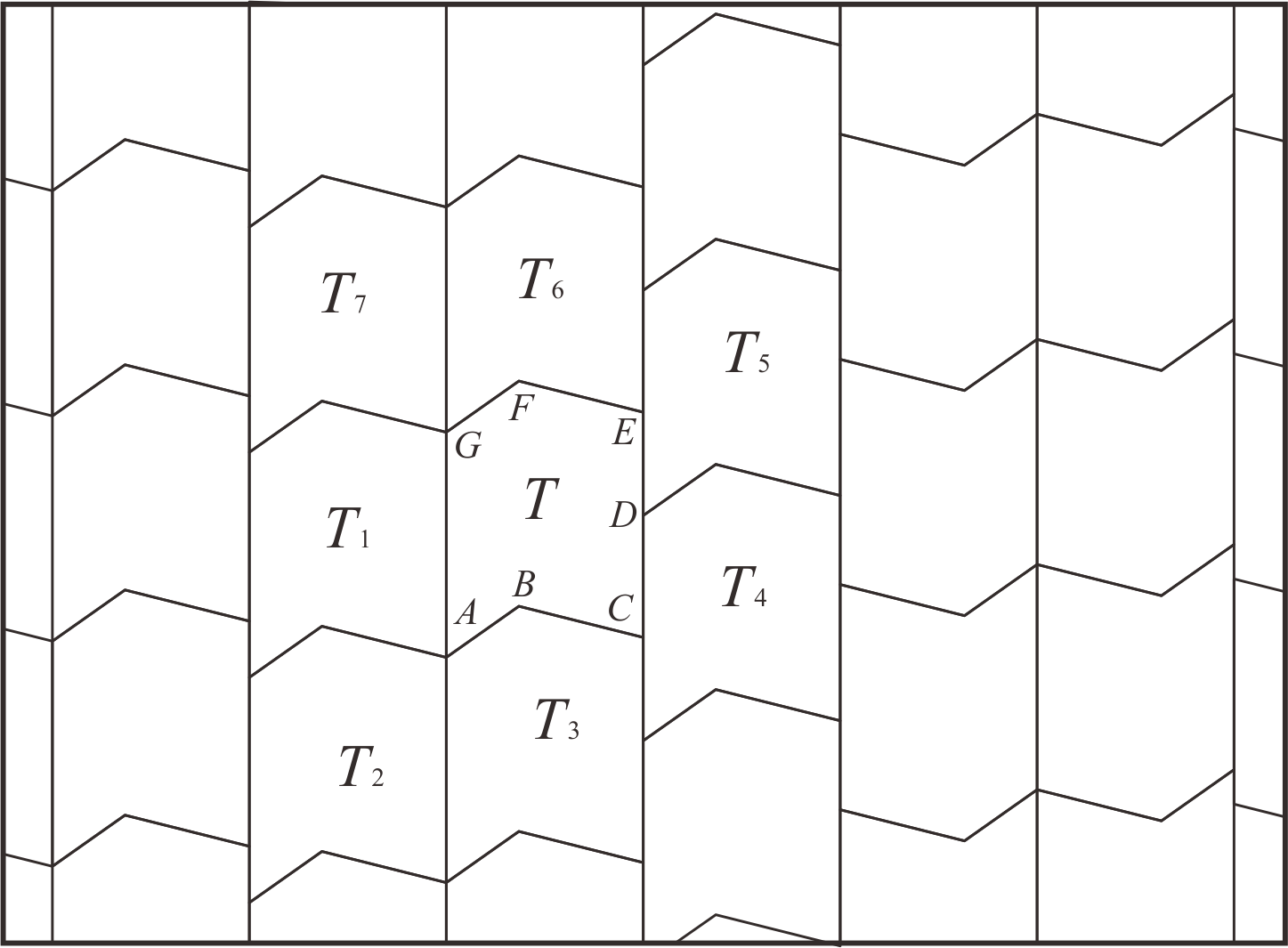} 
  \caption{{\small 
The differences between corners and vertices, sides and edges, 
adjacents, and neighbors. The points $A$, $B$, $C$, $E$, $F$, and $G$ are corners 
of the tile $T$; but $A$, $C$, $D$, $E$, and $G$ are vertices of the tiling (we note 
that the valence of vertices $A$ and $G$ is four, and the valence of vertices 
$C$, $D$, and $E$ is three). The line segments \textit{AB}, \textit{BC}, \textit{CE}, 
\textit{EF}, \textit{FG}, and \textit{GA} are sides of $T$, while \textit{AC}, 
\textit{CD}, \textit{DE}, \textit{EG}, and \textit{GA} are edges of the tiling. 
The tiles $T_{1}$, $T_{3}$, $T_{4}$, $T_{5}$, and $T_{6}$ are adjacents (and 
neighbors) of $T$, whereas tiles $T_{2}$ and $T_{7}$ are neighbors (but not 
adjacents) of $T$ \cite{G_and_S_1987}.
}
\label{fig2}
}
\end{figure}

\section{Preparation}
\label{section2}

Definitions and terms of this section quote from \cite{G_and_S_1987}. 

Terms ``vertices'' and ``edges'' are used by both polygons and tilings. In order 
not to cause confusion, \textit{corners} and \textit{sides} are referred to instead 
of vertices and the edges of polygons, respectively. At a vertex of a polygonal tiling, 
corners of two or more polygons meet and the number of polygons meeting at 
the vertex is called the \textit{valence} of the vertex, and is at least three (see 
Figure~\ref{fig2}). Therefore, an edge-to-edge tiling by polygons is such that the 
corners and sides of the polygons in a tiling coincide with the vertices and edges 
of the tiling.

Two tiles are called \textit{adjacent} if they have an edge in common, and then 
each is called an adjacent of the other. On the other hand, two tiles are called 
\textit{neighbors} if their intersection is nonempty (see Figure~\ref{fig2}).

There exist positive numbers $u$ and $U$ such that any tile contains a certain 
disk of radius $u$ and is contained in a certain disk of radius $U$ in which case 
we say the tiles in tiling are \textit{uniformly bounded}.

A tiling $\Im$ is called \textit{normal} if it satisfies following conditions: 
(i) every tiles of $\Im $ is a topological disk; (ii) the intersection of every two 
tiles of $\Im$ is a connected set, that is, it does not consist of two (or 
more) distinct and disjoint parts; (iii) the tiles of $\Im$ are uniformly 
bounded.

Let $D(r,M)$ be a closed circular disk of radius $r$, centered at any point 
$M$ of the plane. Let us place $D(r,M)$ on a tiling, and let $F_{1}$ and $F_{2}$ 
denote the set of tiles contained in $D(r,M)$ and the set of meeting 
boundary of $D(r,M)$ but not contained in $D(r,M)$, respectively. In 
addition, let $F_{3}$ denote the set of tiles surrounded by these in 
$F_{2}$ but not belonging to $F_{2}$. The set $F_1 \cup F_2 \cup F_3 $ of 
tiles is called the patch $A(r,M)$ of tiles generated by $D(r,M)$.

For a given tiling $\Im$, we denote by $v(r,M)$, $e(r,M)$, and $t(r,M)$ the 
numbers of vertices, edges, and tiles in $A(r,M)$, respectively. The tiling 
$\Im$ is called \textit{balanced} if it is normal and satisfies the following 
condition: the limits

\[
v(\Im ) = \mathop {\lim }\limits_{r \to \infty } \frac{v(r,M)}{t(r,M)}\quad 
\mbox{and}\quad e(\Im ) = \mathop {\lim }\limits_{r \to \infty } 
\frac{e(r,M)}{t(r,M)}
\]

\noindent
exist. Note that $v(r,M) - e(r,M) + t(r,M) = 1$ is called Euler's Theorem 
for Planar Maps.

\begin{sta}[Statement 3.3.13 in \cite{G_and_S_1987}] 
\label{sta1}
Every normal periodic tiling is balanced.
\end{sta}

For a given tiling $\Im$, we write $t_h (r,M)$ for the number of tiles with 
$h$ adjacents in $A(r,M)$, and $v_j (r,M)$ for the numbers of $j$-valent 
vertices in $A(r,M)$. Then the tiling $\Im$ is called \textit{strongly balanced} 
if it is normal and satisfies the following condition: all the limits

\[
t_h (\Im ) = \mathop {\lim }\limits_{r \to \infty } \frac{t_h 
(r,M)}{t(r,M)}\quad \mbox{and}\quad v_j (\Im ) = \mathop {\lim }\limits_{r 
\to \infty } \frac{v_j (r,M)}{t(r,M)}
\]

\noindent
exist. Then,

\[
\sum\limits_{h \ge 3} {t_h (\Im ) = 1} \quad \mbox{and}\quad v(\Im ) 
= \sum\limits_{j \ge 3} {v_j (\Im )}
\]

\noindent
hold. Therefore, every strongly balanced tiling is necessarily balanced.

\begin{sta}[Statement 3.4.8 in \cite{G_and_S_1987}]
\label{sta2}
Every normal periodic tiling is strongly balanced.
\end{sta}

\begin{sta}[Statement 3.5.13 in \cite{G_and_S_1987}]
\label{sta3}
For each strongly balanced tiling $\Im $ we have

\[
\frac{1}{\sum\limits_{j \ge 3} {j \cdot w_j (\Im )} } + 
\frac{1}{\sum\limits_{h \ge 3} {h \cdot t_h (\Im )} } = \frac{1}{2}
\]

\noindent
where

\[
w_j (\Im ) = \frac{v_j (\Im )}{v(\Im )}.
\]
\end{sta}

Thus $w_j (\Im )$ can be interpreted as that fraction of the total number of 
vertices in $\Im$ which have valence $j$, and $\sum\limits_{j \ge 3} {j \cdot 
w_j (\Im )} $ is the \textit{average} valence taken over all the vertices. Since 
$\sum\limits_{h \ge 3} {t_h (\Im ) = 1} $ there is a similar interpretation 
of $\sum\limits_{h \ge 3} {h \cdot t_h (\Im )} $: it is the \textit{average} 
number of adjacents of the tiles, taken over all the tiles in $\Im$. 
Since the valence of the vertex is at least three,

\[
\sum\limits_{j \ge 3} {j \cdot w_j (\Im )} \ge 3.
\]

\section{Proof of Theorem~\ref{thm1}}
\label{section3}

Let $\Im _5^{sb} $ a strongly balanced tiling by convex pentagon. From \cite{Sugimoto_2017rcmm2}, 
following propositions is known.

\begin{prop}\label{prop1}
$ 3 \le \sum\limits_{j \ge 3} {j \cdot w_j (\Im _5^{sb} )} \le \frac{10}{3}$.
\end{prop}

From Proposition~\ref{prop1}, if a convex pentagonal tile can generate $\Im _5^{sb} $, 
then the pentagon must be able to form a vertex of valence three. Here, a 
vertex in which two or more corners are concentrated on the side is called a 
\textit{pseudo-vertex}. That is, the vertices $C$, $D$, and $E$ in Figure~\ref{fig2} 
are the pseudo-vertices of valence three.

 Let $G_{0}$ be a convex pentagon that does not in which the sum of the three 
internal angles is equal to $360^ \circ$. If $G_{0}$ is a convex pentagonal tile that 
can generate  $\Im _5^{sb} $, we first consider that $G_{0}$ has only one corner at 
$90^ \circ$. The reason for this is as follows: From Proposition~\ref{prop1}, if $G_{0}$ 
is a convex pentagonal tile that can generate $\Im _5^{sb} $, $G_{0}$ must form a 
3-valent vertex, and it is the 3-valent pseudo-vertex. Given that the total of the 
internal angles of any convex pentagon is $540^ \circ$, when two corners at the 
3-valent pseudo-vertex are different (for example, a pair of corners $A$ and $B$, 
a pair of corners $A$ and $C$, and so on), the convex pentagon must have a 
property in which the sum of three internal angles is equal to $360^ \circ$. 
Therefore, if $G_{0}$ can form a 3-valent pseudo-vertex, we consider that 
$G_{0}$ has only one corner at $90^ \circ$.

Hereafter, it is assumed that the corner $A$ of $G_{0}$ is $90^ \circ$. If $G_{0}$ 
with $A=90^ \circ$ generates $\Im _5^{sb} $, then the corner $A$ belongs to the 
3-valent pseudo-vertex and the corners $B$, $C$, $D$, and $E$ belong to the 
vertices whose valence is four or more. (Of course, the corner $A$ may also 
have the properties of forming other vertices with valences of four or more. 
However, from Proposition~\ref{prop1}, the corner $A$ must form a 3-valent 
pseudo-vertex.) If the corner $A$ belongs to the 3-valent pseudo-vertex and 
the corners $B$, $C$, $D$, and $E$ belong to the 4-valent vertices in 
$\Im _5^{sb} $ by using $G_{0}$ with $A=90^ \circ$, we consider a model in 
which the ratio between the 3-valent pseudo-vertices and the 4-valent 
vertices is $1: 2$. (It is because that, when the corners $A$ of two pieces 
form one 3-valent pseudo-vertex, each corners $B$, $C$, $D$, and $E$ 
should appear two times at two 4-valent vertices~\cite{Sugi_Ogawa_2006}.) 
For tiling of such models, the average valence of the vertices is given as follow:

\[
\sum\limits_{j \ge 3} {j \cdot w_j (\Im _5^{sb} )} = \frac{3 + 2\times 4}{1 
+ 2} = \frac{11}{3}.
\]

\noindent
This result contradicts Proposition~\ref{prop1}. It is obvious that the average 
valence of the vertices will be even larger if the valence of the vertices to 
which the corners $B$, $C$, $D$, and $E$ belong is further increased. 
Therefore, there is no $G_{0}$ as a convex pentagonal tile that can 
generate  $\Im_5^{sb} $. Then, from Statement~\ref{sta2}, there is no 
$G_{0}$ as a convex pentagonal tile that can generate a periodic tiling.

Thus, a convex pentagonal tile that can generate a periodic tiling  
(i.e., a convex pentagonal tile that admits a periodic tiling) has a 
property in which the sum of three internal angles of the pentagon is equal 
to $360^ \circ$.
\hspace{10cm} $\square$

\section{Conclusions}

We know for a fact that the convex pentagonal tiles belonging to families of 
Types 1--15 admit at least one periodic tiling. On the other hand, there is 
no assurance that all convex pentagonal tiles admit at least one periodic tiling. 
In the solution to convex pentagonal tiling, it is necessary to consider whether 
there is a convex polygonal tile that admits infinitely many tilings of the plane, 
none of which is periodic~\footnote{  A set of prototiles is called aperiodic if 
congruent copies of the prototiles admit infinitely many tilings of the plane, 
none of which are periodic. A tiling that has no periodicity is called nonperiodic. 
On the other hand, a tiling by aperiodic prototiles is called aperiodic. Note 
that, although an aperiodic tiling is a nonperiodic tiling, a nonperiodic tiling 
is not necessarily an aperiodic tiling ~\cite{G_and_S_1987, Hallard_1991, 
Sugimoto_2017rcmm1, wiki_aperiodic_tiling, wiki_aperiodic_set} }
  (i.e., whether there is a convex pentagonal tile which is an aperiodic 
prototile)~\cite{Sugimoto_2017rcmm1, Sugimoto_2017rcmm2}. 
Currently, if a convex pentagonal tile that is an aperiodic prototile exists, 
it is unknown whether the pentagon must have a property in which the 
sum of three internal angles is equal to $360^ \circ$..

\bigskip
\noindent
\textbf{Acknowledgments.} 
The authors would like to thank Professor Yoshio Agaoka of Hiroshima 
University for providing valuable comments.

\end{document}